\documentclass[11pt]{article}
\usepackage{amsmath}
\newtheorem{thm}{Theorem}

\def\eqn#1{(\ref{eq:#1})}

\begin{document}

\makeatletter
\@addtoreset{equation}{section}
\makeatother

\begin{titlepage}

\begin{center}

{\Huge{\bf A double bounded key identity }}\\
{\Huge{\bf for G\"ollnitz's (BIG) partition theorem}}\\

\vspace{2cm}

\large{{\bf Krishnaswami Alladi}} and \large{{\bf Alexander Berkovich}}
\footnote[1]{Research supported in part by a grant from the Number Theory 
Foundation and University of Florida CLAS Research Award}\\
University of Florida, Gainesville, FL 32611 \\
e-mail: alladi@math.ufl.edu; alexb@math.ufl.edu

\end{center}

\vspace{2.5cm}

\begin{abstract}

Given integers $i,j,k,L,M$, we establish a new double bounded $q-$series
identity from which the three parameter $(i,j,k)$ {\it{key identity}} of
Alladi-Andrews-Gordon for G\"ollnitz's (big) theorem follows if
$L,M\to\infty$.  When $L=M$, the identity yields a strong refinement of
G\"ollnitz's theorem with a bound on the parts given by $L$.  This is
the first time a bounded version of G\"ollnitz's (big) theorem has
been proved. This leads to new bounded versions of Jacobi's triple
product identity for theta functions and other fundamental identities.

\end{abstract}

\end{titlepage}

\section{Introduction}

Our goal here is to prove the following double bounded {\it{key identity}}
for G\"ollnitz's (big) partition theorem: {\it{If}} $i,j,k,L,M$, {\it{are
integers, then}}
$$
\sum_{\substack{i,j,k\\ constraints}} q^{T_t+T_{ab}+T_{ac}+T_{bc}}
\left[\begin{array}{c} L-t+a\\ a\end{array}\right]
\left[\begin{array}{c} L-t+b\\ b\end{array}\right]
\left[\begin{array}{c} M-t+c\\ c\end{array}\right]\times
\nonumber
$$
$$
\left[\begin{array}{c} L-t\\ab\end{array}\right]
\left[\begin{array}{c} M-t\\ ac\end{array}\right]
\left[\begin{array}{c} M-t\\ bc\end{array}\right]+
$$
$$
\sum_{\substack{i,j,k\\ constraints}} q^{T_t+T_{ab}+T_{ac}+T_{bc-1}}
\left[\begin{array}{c} L-t+a-1\\ a-1\end{array}\right]
\left[\begin{array}{c} L-t+b\\ b\end{array}\right]
\left[\begin{array}{c} M-t+c\\ c\end{array}\right]\times
$$
$$
\left[\begin{array}{c} L-t\\ ab\end{array}\right]
\left[\begin{array}{c} M-t\\ ac\end{array}\right]
\left[\begin{array}{c} M-t\\ bc-1\end{array}\right]
$$
\begin{equation}
=\sum_{s\ge0}
q^{s(M+2)-T_s+T_{i-s}+T_{j-s}+T_{k-s}}\left[\begin{array}{c} L-s\\ s,
i-s,j-s\end{array}\right] \left[\begin{array}{c}
M-i-j\\k-s\end{array}\right],
\label{eq:1.1}
\end{equation}

\noindent
{\it{where}} $t=a+b+c+ab+ac+bc$ {\it{and the summations on the left in
\eqn{1.1} are over parameters}} $a,b,c,ab,ac$ {\it{and}} $bc$, 
{\it{satisfying the}} $i,j,k$ {\it{constraints}}
\begin{eqnarray}
\left. \begin{array}{l}
 i=a+ab+ac,\\
j=b+ab+bc,\\
k=c+ac+bc.\end{array}\right\}
\label{eq:1.2}
\end{eqnarray}

\noindent
We emphasize that in \eqn{1.1} and everywhere, $ab$ is a
variable, and is {\underbar{not equal to $a$ times $b$}}, with similar
interpretation for $ac$ and $bc$. The role of the variables $a,b,\dots,bc$
will become clear in the sequel.

In \eqn{1.1} and in what follows, $T_n=n(n+1)/2$, and the $q-$binomial and
$q-$multinomial coefficients are defined by
\begin{equation}
\left[\begin{array}{c} n+m\\ n\end{array}\right]=\left[\begin{array}{c}
n+m\\ n\end{array}\right]_q= \left\{\begin{array}{ll}
\frac{(q^{m+1})_n}{(q)_n}, & \mbox{if } n\ge0,\\ 0, & \mbox{otherwise},
\end{array}\right.
\label{eq:1.3}
\end{equation}
and
\begin{eqnarray}
\left[\begin{array}{c} L\\
a,b,c,\dots\end{array}\right]&=&\left[\begin{array}{c} L\\
a\end{array}\right]\left[\begin{array}{c} L-a\\ b\end{array}\right]
\left[\begin{array}{c} L-a-b\\ c\end{array}\right]\dots \nonumber\\
&=& \left\{\begin{array}{ll}
\frac{(q^{1+L-a-b-c\dots})_{a+b+c\dots}}{(q)_a(q)_b(q)_c\dots}, &
\mbox{if } a\ge0,b\ge0,c\ge0,\dots,\\ 0, &
\mbox{otherwise},\end{array}\right.
\label{eq:1.4}
\end{eqnarray}
where the symbols $(a)_n$ are given by
\begin{equation}
(a;q)_n=(a)_n=\left\{\begin{array}{lll}
\prod^{n-1}_{j=0}(1-aq^j), & \mbox{ if } & n>0,\\ 1, & \mbox{ if } &
n=0,\\ \prod^{-n}_{j=1} (1-aq^{-j})^{-1}, & \mbox{ if } &
n<0.\end{array} \right.
\label{eq:1.5}
\end{equation}

The connections between \eqn{1.1} and the partition theorem of G\"ollnitz
\cite{G} will be explained subsequently.  Note that when $L,M\to \infty$, only 
the term corresponding to $s=0$ on the right hand side of \eqn{1.1} survives,
and so \eqn{1.1} reduces to
\begin{equation}
\sum_{\substack{i,j,k\\ constraints}} \frac{q^{T_t+T_{ab}+T_{ac}+T_{bc-1}}
(1-q^a+q^{a+bc})} {(q)_a(q)_b(q)_c(q)_{ab}(q)_{ac}(q)_{bc}}=
\frac{q^{T_i+T_j+T_k}}{(q)_i(q)_j(q)_k},
\label{eq:1.6}
\end{equation}
which is the three parameter key identity for G\"ollnitz's theorem due
to Alladi-Andrews-Gordon \cite{AAG}. \\
If any one of the parameters $i,j,k$ is set equal to 0, then \eqn{1.1}
reduces to the double bounded key identity for Schur's theorem we have
recently established \cite{AB1}.  For instance, with $i=0$, \eqn{1.1} becomes
\begin{equation}
\sum_{\substack{j=b+bc\\ k=c+bc}} q^{T_{b+c+bc}+T_{bc}}
\left[\begin{array}{c} L-k\\ j-bc\end{array}\right]
\left[\begin{array}{c} M-j\\ k-bc\end{array}\right]
\left[\begin{array}{c} M-b-c-bc\\ bc\end{array}\right]=q^{T_j+T_k}
\left[\begin{array}{c} L\\j\end{array}\right]
\left[\begin{array}{c} M-j\\ k\end{array}\right].
\label{eq:1.7}
\end{equation}

Our proof of \eqn{1.1} has two parts.  Denoting the left hand side of \eqn{1.1}
by $g_{i,j,k}(L,M)$ and the right hand side of \eqn{1.1} by
$p_{i,j,k}(L,M)$, we first show in \S2 that the functions
$g_{i,j,k}(L,M)$ and $p_{i,j,k}(L,M)$ satisfy identical second order
recurrences in $L$.  To complete the proof of the equality
\begin{equation}
g_{i,j,k}(L,M)=p_{i,j,k}(L,M)
\label{eq:1.8}
\end{equation}
we show in \S3 that both functions satisfy the same initial conditions
\begin{equation}
g_{i,j,k}(i+j-1,M)=p_{i,j,k}(i+j-1,M).
\label{eq:1.9}
\end{equation}
This is not as easy as it sounds; the proof of \eqn{1.9} in \S3 requires the
use of Jackson's $q-$analog of Dougall's summation.
When $L=M$, the right hand side of \eqn{1.1} can be evaluated elegantly in
terms of a product of $q-$binomial coefficients with cyclic dependence on
$i,j$, and $k$ (see \S4).  This has a nice partition interpretation
yielding a strong refinement of G\"ollnitz's theorem with a bound on the
size of the parts.  To the best of our knowledge, this is the first time
a bounded version of G\"ollnitz's theorem has been found.  There are a
number of important consequences of this theorem one of which is a new 
finite version of Jacobi's triple product identity which is stated as 
identity \eqn{5.2} in \S5 (also see \eqn{5.3}, \eqn{5.5}); the proof of 
\eqn{5.2} and finite versions of 
many other fundamental results in the theory of partitions and $q$-series 
will be given elsewhere \cite{AB3}, \cite{AB4}.
In \S5 some problems for further investigation motivated by this work
are briefly indicated as well.  Finally, certain technical details
pertaining to recurrences for $q-$multinomial coefficients and to partition
theoretical interpretation of \eqn{1.1} are relegated to Appendix A and B,
respectively.

\section{Recurrences}

Define for integers $i,j,k,\delta,L,M$, the sum
\begin{equation}
X_{i,j,k}(L,M)=\sum_{\substack{i,j,k\\ constraints}} q^{T_t+T_{ab}}
\left[\begin{array}{c} L-t+a-\delta\\ a-\delta\end{array}\right]
\left[\begin{array}{c} L-t+b\\ b\end{array}\right]
\left[\begin{array}{c}
L-t\\ ab\end{array}\right] f(M-t;c,ac,bc),
\label{eq:2.1}
\end{equation}
where $t=a+b+c+ab+ac+bc$ as before, and the explicit form of the
function $f(M;c,ac,bc)$ will not be required for the recurrences. However,
it is important that $f(M;c,ac,bc)$ does not depend on $L;a,b,ab$.
We wish to show that $X_{i,j,k}(L,M)$ satisfies the following second order
recurrence in $L$:
$$
X_{i,j,k}(L,M) = X_{i,j,k}(L-1,M)+q^LX_{i-1,j,k} (L-1,M-1)+
q^LX_{i,j-1,k} (L-1,M-1)
$$
\begin{equation}
+q^LX_{i-1,j-1,k}(L-2,M-1)-q^{2L-1}X_{i-1,j-1,k}(L-2,M-2).
\label{eq:2.2}
\end{equation}
To this end we will use repeatedly the $q-$binomial recurrence
\begin{equation}
\left[\begin{array}{c} n+m\\ n\end{array}\right]=\left[\begin{array}{c}
n+m-1\\ n\end{array}\right] + q^m\left[\begin{array}{c} n-1+m\\
n-1\end{array}\right]
\label{eq:2.3}
\end{equation}
which holds for all integers $m,n$, to expand the right hand side of
\eqn{2.1} in a telescopic fashion as follows:
$$
\sum_{\substack{i,j,k\\ constraints}} q^{T_t+T_{ab}} \left[\begin{array}{c}
L-1-t+a-\delta\\ a-\delta\end{array}\right] \left[\begin{array}{c}
L-1-t+b\\ b\end{array}\right] \left[\begin{array}{c} L-1-t\\
ab\end{array}\right] f(M-t;c,ac,bc)
$$
$$
+\sum_{\substack{i,j,k\\ constraints}} q^{T_{t-1}+T_{ab-1}+L}
\left[\begin{array}{c} L-1-t+a-\delta\\ a-\delta\end{array}\right]
\left[\begin{array}{c} L-1-t+b\\ b\end{array}\right]
\left[\begin{array}{c} L-1-t\\ ab-1\end{array}\right]
f(M-t;c,ac,bc)
$$
$$
+\sum_{\substack{i,j,k\\ constraints}} q^{T_{t-1}+T_{ab}+L}
\left[\begin{array}{c} L-1-t+a-\delta\\ a-\delta\end{array}\right]
\left[\begin{array}{c} L-t+b-1\\ b-1\end{array}\right]
\left[\begin{array}{c} L-t\\ ab\end{array}\right]
f(M-t;c,ac,bc)
$$
\begin{equation}
+\sum_{\substack{i,j,k\\ constraints}} q^{T_{t-1}+T_{ab}+L}
\left[\begin{array}{c}
L-t+a-1-\delta\\ a-1-\delta\end{array}\right] \left[\begin{array}{c}
L-t+b\\ b\end{array}\right] \left[\begin{array}{c} L-t\\
ab\end{array}\right] f(M-t;c,ac,bc).
\label{eq:2.4}
\end{equation}

\noindent
Let us denote each of the four sums in \eqn{2.4} by $\Sigma_1,\Sigma_2,\Sigma_3$,
and $\Sigma_4$, respectively.  To see that \eqn{2.4} is an expansion of \eqn{2.1},
we merge $\Sigma_1$ and $\Sigma_2$ in \eqn{2.4} into a single sum with the aid
of \eqn{2.3}.  This single sum can in turn be merged with $\Sigma_3$ in \eqn{2.4}
using \eqn{2.3}, and this finally can be merged with $\Sigma_4$ in \eqn{2.4} to
yield \eqn{2.1}.  (This telescopic expansion technique was introduced in \cite{BER}
and later was used
extensively by Berkovich-McCoy-Schilling \cite{BMS} and by Schilling-Warnaar
\cite{SW}.)
From the definition of $X_{i,j,k}(L,M)$ it is clear that
\begin{equation}
\Sigma_1=X_{i,j,k}(L-1,M).
\label{eq:2.5}
\end{equation}
If we perform the change $ab\mapsto ab+1$ in $\Sigma_2$, then $t\mapsto
t+1$, $i\mapsto i-1$, $j\mapsto j-1$, and so
$$
\Sigma_2 = \sum_{\substack{i-1,j-1,k\\ constraints}} q^{T_t+T_{ab}+L}
\left[\begin{array}{c} L-2-t+a-\delta\\ a-\delta\end{array}\right]
\left[\begin{array}{c} L-2-t+b\\ b\end{array}\right] \left[\begin{array}{c}
L-2-t\\ ab\end{array}\right] \times
$$
\begin{equation}
f(M-1-t;c,ac,bc)=q^LX_{i-1,j-1,k}(L-2,M-1).
\label{eq:2.6}
\end{equation}
Similarly, replacing $a$ by $a+1$ in $\sum_4$, we obtain
\begin{equation}
\Sigma_4=q^LX_{i-1,j,k}(L-1,M-1).
\label{eq:2.7}
\end{equation}

With regard to $\Sigma_3$, we write it as a difference to recognize it as
$$
\Sigma_3=\sum_{\substack{i,j,k\\ constraints}} q^{T_{t-1}+T_{ab}+L}
\left[\begin{array}{c} L-t+a-\delta\\ a-\delta\end{array}\right]
\left[\begin{array}{c} L-t+b-1\\ b-1\end{array}\right]
\left[\begin{array}{c} L-t\\ ab\end{array}\right]
f(M-t;c,ac,bc)
$$
$$
-\sum_{\substack{i,j,k\\ constraints}} q^{T_{t-2}+T_{ab}+2L-1}
\left[\begin{array}{c} L-t+a-1-\delta\\ a-1-\delta\end{array}\right]
\left[\begin{array}{c} L-t+b-1\\ b-1\end{array}\right]
\left[\begin{array}{c} L-t\\ ab\end{array}\right]
f(M-t;c,ac,bc)
$$
\begin{equation}
=q^LX_{i,j-1,k}(L-1,M-1)-q^{2L-1}X_{i-1,j-1,k}(L-2,M-2).
\label{eq:2.8}
\end{equation}

\noindent
Observe that $g_{i,j,k}(L,M)$, the left hand side of \eqn{1.1}, is a
sum of two functions $X_{i,j,k}(L,M)$, one with $\delta=0$, and
the other with $\delta=1$, and with $f(M-t;c,ac,bc)$ suitably
identified.  So it follows that
$$
g_{i,j,k}(L,M)=g_{i,j,k}(L-1,M)+q^Lg_{i-1,j,k}(L-1,M-1)+q^Lg_{i,j-1,k}
(L-1,M-1)
$$
\begin{equation}
+q^Lg_{i-1,j-1,k} (L-2,M-1)-q^{2L-1}g_{i-1,j-1,k} (L-2,M-2).
\label{eq:2.9}
\end{equation}
In \cite{A1}, Andrews had derived a fourth order recursion relation in
$L$ for $g_{i,j,k}(L,L)$. His recurrence can be generalized as
$$
g_{i,j,k}(L,M)=g_{i,j,k}(L-1,M-1)+(q^L g_{i-1,j,k}(L-1,M-1)+
q^L g_{i,j-1,k}(L-1,M-1)+ q^M g_{i,j,k-1}(L-1,M-1))
$$
$$
+(1-q^{L-1})(q^L g_{i-1,j-1,k}(L-2,M-2)+q^M g_{i-1,j,k-1}(L-2,M-2)
+q^M g_{i,j-1,k-1}(L-2,M-2))
$$
$$
+q^{2L+M-3}g_{i-1,j-1,k-1}(L-3,M-3)
$$
$$
+(q^{L+M-1}g_{i-2,j-1,k-1}(L-3,M-3)+q^{L+M-1} g_{i-1,j-2,k-1}(L-3,M-3)
+q^{2M-1}g_{i-1,j-1,k-2}(L-3,M-3))
$$
$$
+q^{L+2M-3}g_{i-2,j-2,k-2}(L-4,M-4).
\nonumber
$$
The introduction of an extra parameter $M$
has enabled us to bring down the order of the recursion relation in $L$
for $g_{i,j,k}(L,M)$ to just two.  

Next, we claim that $p_{i,j,k}(L,M)$, the right hand side of \eqn{1.1},
satisfies the same recurrence, namely,
$$
p_{i,j,k}(L,M) = p_{i,j,k}(L-1,M)+q^Lp_{i-1,j,k}
(L-1,M-1)+q^Lp_{i,j-1,k} (L-1,M-1)+
$$
\begin{equation}
+ q^Lp_{i-1,j-1,k} (L-2,M-1)-q^{2L-1} p_{i-1,j-1,k} (L-2,M-2).
\label{eq:2.10}
\end{equation}
For this purpose we employ the following recursion relation for the
$q-$multinomial coefficients (see Appendix A for a proof):
$$
\left[\begin{array}{c} L-s\\ s,i-s,j-s\end{array}\right]
=\left[\begin{array}{c} L-1-s\\ s,i-s, j-s\end{array}\right] +
q^{L-i}\left[\begin{array}{c} L-1-s\\ s, i-1-s,j-s\end{array}\right]+
$$
$$
q^{L-j}\left[\begin{array}{c} L-1-s\\ s,i-s,j-1-s\end{array}\right]+
q^{L-i-j}\left[\begin{array}{c} L-1-s\\ s-1,i-s,j-s\end{array}\right]+
$$
\begin{equation}
q^{L+s-i-j} \left[\begin{array}{c} L-2-s\\ s,i-1-s,j-1-s\end{array}\right]
- q^{2L-1-i-j} \left[\begin{array}{c}L-2-s\\
s,i-1-s,j-1-s\end{array}\right].
\label{eq:2.11}
\end{equation}
Substituting \eqn{2.11} into the right hand side of \eqn{1.1}, we see that
$$
p_{i,j,k}(L,M)=p_{i,j,k}(L-1,M)+q^Lp_{i-1,j,k}
(L-1,M-1)+q^Lp_{i,j-1,k} (L-1,M-1)
$$
$$
\sum_{s\ge0} q^{L-i-j+s(M+2)-T_s+T_{i-s}+T_{j-s}+T_{k-s}}
\left[\begin{array}{c} L-1-s\\s-1,i-s,j-s\end{array}\right]
\left[\begin{array}{c} M-i-j\\k-s\end{array}\right]+
$$
$$
\sum_{s\ge0} q^{L+s-i-j+s(M+2)-T_s+T_{i-s}+T_{j-s}+T_{k-s}}
\left[\begin{array}{c} L-2-s\\s,i-1-s,j-1-s\end{array}\right]
\left[\begin{array}{c} M-i-j\\k-s\end{array}\right]
$$
\begin{equation}
-q^{2L-1}p_{i-1,j-1,k} (L-2,M-2).
\label{eq:2.12}
\end{equation}
Now we replace $s$ by $s+1$ in the first sum in \eqn{2.12}.  This enables us
to merge this sum with the second sum in \eqn{2.12} to obtain
$$
q^L \sum_{s\ge0} q^{s(M-1+2)-T_s+T_{i-1-s}+T_{j-1-s}+T_{k-s}}
\left[\begin{array}{c} L-2-s\\s,i-1-s,j-1-s\end{array}\right]
\left[\begin{array}{c} M-1-(i-1)-(j-1)\\k-s\end{array}\right]
$$
\begin{equation}
=q^Lp_{i-1,j-1,k}(L-2,M-1).
\label{eq:2.13}
\end{equation}
The recurrence \eqn{2.10} follows from \eqn{2.12} and \eqn{2.13}.

\section{The boundary identity}

Having established that $g_{i,j,k}(L,M)$ and $p_{i,j,k}(L,M)$ satisfy
identical recurrences \eqn{2.9} and \eqn{2.10}, we note now that
\begin{equation}
g_{i,j,k}(L,M)=p_{i,j,k}(L,M)=0
\label{eq:3.1}
\end{equation}
if any one of the parameters $i,j,k$ is negative.  Thus if we show that
the boundary identity \eqn{1.9} is true, then we can conclude that
\begin{equation}
g_{i,j,k}(L,M)=p_{i,j,k}(L,M), \mbox{ }
\forall(i,j,k,L,M)\in\mathbf{Z}^5.
\label{eq:3.2}
\end{equation}

A few comments are in order concerning the nonstandard choice of the
diagonal boundary $L=i+j-1$.  The conventional choice $L=0,1$ leads to
difficulties because the terms of \eqn{1.1} do not collapse in these cases.
Moreover, the truth of \eqn{1.1} for $L=0,1$, leads us to conclude its
validity only for $L\ge0$. Consequently, the case $L<0$ of \eqn{1.1} which is
highly nontrivial would not be covered.  On the other hand the choice
$L=i+j-1$ enables us to prove \eqn{1.1} for all $L \in \mathbf{Z}$.
To provide additional motivation for the choice $L=i+j-1$, we now show
that $p_{i,j,k}(L,M)$ collapses radically in this case.  Indeed,
$$
p_{i,j,k}(i+j-1,M) =
\sum_{s\ge0}q^{s(M+2)-T_s+T_{i-s}+T_{j-s}+T_{k-s}}
\left[\begin{array}{c} i+j-s-1\\s\end{array}\right]
\left[\begin{array}{c} i-s+j-s-1\\i-s\end{array}\right]\times
$$
\begin{equation}
\left[\begin{array}{c} j-s-1\\j-s\end{array}\right]
\left[\begin{array}{c} M-i-j\\k-s\end{array}\right]
= \delta_{i,0}\delta_{j,0} q^{T_k}
\left[\begin{array}{c} M-i-j\\k\end{array}\right]
\label{eq:3.3}
\end{equation}
by noticing that
\begin{equation}
\left[\begin{array}{c} j-s-1\\ j-s\end{array}\right]=\delta_{j,s};
\left[\begin{array}{c} i-j-1\\ i-j\end{array}\right]=\delta_{i,j};
\left[\begin{array}{c} i-1\\ i\end{array}\right]=\delta_{i,0}
\label{eq:3.4}
\end{equation}
where
$$
\delta_{i,j}=\left\{\begin{array}{ll} 1, & \mbox{if } i=j,\\ 0, &
\mbox{otherwise}.\end{array}\right.
$$
Thus the boundary identity \eqn{1.9} can be stated as
\begin{equation}
g_{i,j,k}(i+j-1,M)=\delta_{i,0}\delta_{j,0}q^{T_k}
\left[\begin{array}{c}\Delta\\ k\end{array}\right],
\label{eq:3.5}
\end{equation}
where $\Delta=M-i-j$.
Next, by repeated use of the $q-$binomial formula
\begin{equation}
\left[\begin{array}{c} -\alpha\\k\end{array}\right]=(-1)^k
\left[\begin{array}{c} k+\alpha-1\\k\end{array}\right] q^{-\alpha k-T_{k-1}}
\label{eq:3.6}
\end{equation}
(see Gasper and Rahman \cite{GR}, formula (I.44)), we may rewrite \eqn{3.5} as
$$
\sum_{\substack{i,j,k\\ constraints}} (-1)^{a+b+ab}
q^{T_\tau+T_{\Gamma}+T_{ac}+T_{bc}+T_{a-1} +T_{b-1}
+\Gamma(ab+ac+bc)}\times
$$
$$
\qquad\qquad \left[\begin{array}{c}\Gamma\\ a \end{array}\right]
\left[\begin{array}{c} \Gamma\\ b \end{array}\right]
\left[\begin{array}{c} \Delta\\ \Gamma\end{array}\right]
\left[\begin{array}{c} ab+\Delta\\ ab\end{array}\right]
\left[\begin{array}{c}\Delta-\Gamma\\ ac \end{array}\right]
\left[\begin{array}{c} \Delta-\Gamma\\ bc\end{array}\right]+
$$
$$
\sum_{\substack{i,j,k\\ constraints}} (-1)^{a-1+b+ab}
q^{T_\tau+T_\Gamma+T_{ac}+T_{bc-1}+T_{a-2} +T_{b-1}
+\Gamma(1+ab+ac+bc)}\times
$$
$$
\qquad\qquad \left[\begin{array}{c}\Gamma\\ a-1 \end{array}\right]
\left[\begin{array}{c} \Gamma\\ b \end{array}\right]
\left[\begin{array}{c} \Delta\\ \Gamma\end{array}\right]
\left[\begin{array}{c} ab+\Delta\\ ab\end{array}\right]
\left[\begin{array}{c}\Delta-\Gamma\\ ac \end{array}\right]
\left[\begin{array}{c} \Delta-\Gamma\\ bc-1\end{array}\right]
$$
\begin{equation}
\qquad\qquad\qquad\qquad\qquad=\delta_{i,0}\delta_{j,0}q^{T_k}
\left[\begin{array}{c} \Delta\\k\end{array}\right],
\label{eq:3.7}
\end{equation}

\noindent
where
\begin{equation}
\Gamma=c-ab\mbox{ and } \tau=a+b+2ab+ac+bc.
\label{eq:3.8}
\end{equation}

It is convenient to treat $\Delta$ as an independent parameter and $\Gamma$ as
an independent summation variable. In this case the constraints in \eqn{3.7}
become \begin{equation}
\left.\begin{array}{l}
i=a+ab+ac,\\ j=b+ab+ac,\\ k=\Gamma+ab+ac+bc.\end{array}\right\}
\label{eq:3.9}
\end{equation}

\noindent
Next, multiply both sides of \eqn{3.7} by $A^iB^jC^k$ and sum over $i,j,k$.
For the right hand side we get immediately
\begin{equation}
\sum_{i,j,k\ge0}\delta_{i,0} \delta_{j,0} q^{T_k}\left[\begin{array}{c}
\Delta\\ k\end{array}\right] A^i B^jC^k=(-Cq)_\Delta.
\label{eq:3.10}
\end{equation}

To treat the left hand side of \eqn{3.7}, we get rid of the condition on
$\tau$ in \eqn{3.8} and rewrite it as
$$
[\omega^0]
\left\{\theta(\omega,q) \sum_{a,b,\Gamma,ab,ac,bc}
q^{T_\Gamma+T_{ac}+T_{b-1}}\right.
\left(-\frac A\omega\right)^a
\left(-\frac B\omega\right)^b
\left(-\frac{ABC q^\Gamma}{\omega^2}\right)^{ab}
C^\Gamma
\left(\frac{ACq^\Gamma}\omega\right)^{ac}
\left(\frac{BCq^\Gamma}\omega\right)^{bc}\times
$$
$$
\left[\begin{array}{c} \Gamma\\ b \end{array}\right]
\left[\begin{array}{c} \Delta\\ \Gamma \end{array}\right]
\left[\begin{array}{c} \Delta+ab\\ ab\end{array}\right]
\left[\begin{array}{c} \Delta-\Gamma\\ ac\end{array}\right]\times
$$
\begin{equation}
\left.\left(q^{T_{bc}+T_{a-1}}
\left[\begin{array}{c} \Gamma\\ a\end{array}\right]
\left[\begin{array}{c} \Delta-\Gamma\\ bc\end{array}\right]
-q^{\Gamma+T_{bc-1}+T_{a-2}}
\left[\begin{array}{c} \Gamma\\ a-1\end{array}\right]
\left[\begin{array}{c} \Delta-\Gamma\\
bc-1\end{array}\right]\right)\right\},
\label{eq:3.11}
\end{equation}

\noindent
where
$$
\theta(\omega,q)=\sum^\infty_{\tau=-\infty}\omega^\tau q^{T_\tau},
$$
and $[\omega^m]f(\omega)$ is the coefficient of $\omega^m$ in the
Laurent expansion of $f(\omega)$.
Thanks to the two $q-$binomial theorems
\begin{equation}
\sum_{n\ge0} z^nq^{T_n}\left[\begin{array}{c} \Delta\\
n\end{array}\right]=(-zq)_\Delta,
\label{eq:3.12}
\end{equation}
and
\begin{equation}
\sum_{n\ge0} z^n \left[\begin{array}{c} \Delta+n\\
n\end{array}\right]=\frac1{(z)_{\Delta+1}},
\label{eq:3.13}
\end{equation}
we can evaluate the summations in \eqn{3.11} over the variables $a,b,ab,ac$,
and $bc$, to cast the left hand side of \eqn{3.7} as
\begin{equation}
[\omega^0]\left\{\theta(\omega,q) \sum_{\Gamma\ge0}
\frac{(1+\frac{ABC}{\omega^2}q^{2\Gamma})}{\left(-\frac{ABC}{\omega^2}
q^\Gamma\right)_{\Delta+1}}
\left[\begin{array}{c} \Delta\\ \Gamma\end{array}\right] C^\Gamma
q^{T_\Gamma}\left(\frac {A}{\omega}\right)_\Gamma
\left(\frac {B}{\omega}\right)_\Gamma
\left(-\frac {AC}{\omega}q^{1+\Gamma}\right)_{\Delta- \Gamma}
\left(-\frac{BC}{\omega}q^{1+\Gamma}\right)_{\Delta-\Gamma}\right\}.
\label{eq:3.14}
\end{equation}
We would like to write \eqn{3.14} in $q-$hypergeometric form.  This can be
done with the aid of the following formulas:
\begin{eqnarray}
&(i) & \left[\begin{array}{c} \Delta\\ \Gamma\end{array}\right]=
\frac{(q^{-\Delta})_\Gamma} {(q)_\Gamma}
(-q^\Delta)^\Gamma q^{-T_{\Gamma-1}},\nonumber\\
&(ii) &
(xq^\Gamma)_{\Delta-\Gamma}=\frac{(x)_\Delta}{(x)_\Gamma},\nonumber\\
&(iii) &
\frac{1}{(xq^\Gamma)_\Delta}=\frac{(x)_\Gamma}{(x)_\Delta(xq^\Delta)_\Gamma},
\nonumber\\
&(iv) & 1+\frac{ABC}{\omega^2}q^{2\Gamma}=
\left(1+\frac{ABC}{\omega^2}\right)
\frac{\left(q\sqrt{\frac{-ABC}{\omega^2}},
-q\sqrt{\frac{-ABC}{\omega^2}}\right)_\Gamma}
{\left(\sqrt{\frac{-ABC}{\omega^2}},
-\sqrt{\frac{-ABC}{\omega^2}}\right)_\Gamma},
\label{eq:3.15}
\end{eqnarray}
where
\begin{equation}
(a_1,a_2,\dots,a_r;q)_m=(a_1,a_2,\dots,a_r)_m=
(a_1)_m(a_2)_m\dots(a_r)_m.
\label{eq:3.16}
\end{equation}
Thus the expression in \eqn{3.14} is
\begin{equation}
[\omega^0]\left\{\theta(\omega,q)
\frac{\left(\frac{-ACq}\omega,\frac{-BCq}\omega\right)_\Delta}{\left(\frac{-
ABC}{\omega^2}
q\right)_\Delta} \mbox{ }_6\phi_5 \left(\begin{array}{l} y, q\sqrt y,-q\sqrt
y, \frac A\omega,\frac B\omega,q^{-\Delta}\\
\sqrt y,-\sqrt y, -\frac{AC}{\omega} q, -\frac{BC}\omega q, yq^{\Delta+1}
\end{array}; q,-cq^{\Delta+1}\right)\right\},
\label{eq:3.17}
\end{equation}
where $y=-\frac{ABC}{\omega^2}$, and we have made use of standard notation
\begin{equation}
{}_{r+1}\phi_r\left(\begin{array}{l} a_1,a_2,\dots,a_{r+1}\\ b_1,b_2,\dots,
b_r\end{array}; q,z\right)= \sum_{n\ge0}
\frac{(a_1,a_2,\dots,a_{r+1})_n} {(q,b_1,b_2,\dots,b_r)_n} z^n.
\label{eq:3.18}
\end{equation}
Actually the ${}_6\phi_5$ in \eqn{3.17} can be evaluated by Jackson's
$q-$analog of Dougall's summation (see \cite{GR}, formula (II.21)) to be
\begin{equation}
\frac{\left(-\frac{ABC}{\omega^2} q, -
Cq\right)_\Delta}{\left(-\frac{AC}\omega q, -\frac{BC}\omega
q\right)_\Delta}.
\label{eq:3.19}
\end{equation}
Finally, combining \eqn{3.7}, \eqn{3.10}, \eqn{3.11}, \eqn{3.14},\eqn{3.17} and \eqn{3.19}, we
can rewrite \eqn{3.5} as
$$
[\omega^0](\theta(\omega,q)\cdot (-C q)_\Delta)=(-Cq)_\Delta,
$$
which is obviously true because
$$
[\omega^0]\theta(\omega,q)=1.
$$
Thus, we have completed the proof of the boundary identity \eqn{1.9} and
consequently the truth of \eqn{1.8} (and \eqn{1.1}) is established.

\section{A bounded version of G\"ollnitz's partition theorem}

In this section we assume that $L=M$, and for this case we first
establish the representation
\begin{equation}
p_{i,j,k}(L,L)=q^{T_i+T_j+T_k}\left[\begin{array}{c} L-k\\
i\end{array}\right] \left[\begin{array}{c} L-i\\ j\end{array}\right]
\left[\begin{array}{c} L-j\\ k\end{array}\right].
\label{eq:4.1}
\end{equation}
We then discuss the partition interpretation of the identity
\begin{equation}
g_{i,j,k}(L,L)=p_{i,j,k}(L,L).
\label{eq:4.2}
\end{equation}

First note that
\begin{equation}
\left[\begin{array}{c} L-s\\ s,i-s,j-s\end{array}\right]
\left[\begin{array}{c} L-i-j\\ k-s\end{array}\right]=
\left[\begin{array}{c} L-s\\ s,i-s,j-s,k-s\end{array}\right].
\label{eq:4.3}
\end{equation}
With \eqn{4.3} in mind, we rewrite $p_{i,j,k}(L,L)$ as
$$
p_{i,j,k}(L,L) = \underset{\ell\to L}{\lim}\sum_{s\ge0}
q^{s(\ell+2)-T_s+T_{i-s}+T_{j-s}+T_{k-s}} \left[\begin{array}{c} \ell-s\\
s,i-s,j-s,k-s\end{array}\right]=\nonumber\\
$$
\begin{eqnarray}
\underset{\ell\to L}{\lim} \frac{q^{T_i+T_j+T_k}}{(q)_i(q)_j(q)_k}
(q^{1+\ell-i-j-k})_{i+j+k}\times{ }_3\phi_2 \left(\begin{array}{l}
q^{-i},q^{-j},q^{-k}\\ q^{1+\ell-i-j-k}, q^{-\ell}\end{array};
q,q\right),
\label{eq:4.4}
\end{eqnarray}
where we have used the limit definition to make sure that all objects in
\eqn{4.4} are well defined.  It turns out that by the use of the
$q-$Pfaff-Saalsch\"utz summation (see Gasper and Rahman \cite{GR}, eqn.(II.12))
\begin{equation}
{}_3\phi_2 \left(\begin{array}{l} q^{-i},q^{-j},q^{-k}\\
q^{1+\ell-i-j-k},q^{-\ell}\end{array}; q,q\right)=
\frac{(q^{1+\ell-j-k}, q^{1+\ell-i-k})_k}{(q^{1+\ell-i-j-k},
q^{1+\ell-k})_k},
\label{eq:4.5}
\end{equation}

\noindent
and so \eqn{4.4} becomes
\begin{equation}
p_{i,j,k}(L,L)=\lim_{\ell\to L} \frac{q^{T_i+T_j+T_k}}{(q)_i(q)_j}
\frac{(q^{1+\ell-i-j-k})_{i+j+k} (q^{1+\ell-i-k})_k}
{(q^{1+\ell-i-j-k})_k(q^{1+\ell-k})_k} \left[\begin{array}{c} \ell-j\\
k\end{array}\right].
\label{eq:4.6}
\end{equation}
Finally it can be shown by repeated use of the formula (ii) of \eqn{3.15}
that \eqn{4.6} yields (we omit the lengthy details of this calculation)
$$
p_{i,j,k}(L,L)=\lim_{\ell\to L}q^{T_i+T_j+T_k}
\left[\begin{array}{c} \ell-k\\i\end{array}\right]
\left[\begin{array}{c} \ell-i\\ j\end{array}\right]
\left[\begin{array}{c} \ell-j\\ k\end{array}\right],
$$
which is \eqn{4.1}, thus completing the proof.

{\bf{Remark 1}}: Zeilberger \cite{Z} has given a combinatorial proof of the 
q-Pfaff-Saalsch\"utz summation. His parameters translate to ours in (4.3) 
by suitable change of variables. 

Now when
\begin{equation}
L\ge\max(i+j, j+k, k+i),
\label{eq:4.7}
\end{equation}
$p_{i,j,k}(L,L)$ given by \eqn{4.1} can be interpreted as the generating
function of partitions $\pi$ whose parts occur in three (primary) colors
$\mathbf{A},\mathbf{B},\mathbf{C}$ ordered as
\begin{equation}
\mathbf{A}<\mathbf{B}<\mathbf{C}
\label{eq:4.8}
\end{equation}
such that parts in the same color are distinct and
\begin{equation}
\left.\begin{array}{llll}
\nu(\mathbf{A};\pi)= & \nu(\mathbf{A})=i, & \lambda(\mathbf{A},\pi)= &
\lambda(\mathbf{A})\le L-k,\\
& \nu(\mathbf{B})=j, & &  \lambda(\mathbf{B})\le L-i,\\
& \nu(\mathbf{C})=k, & & \lambda(\mathbf{C})\le L-j,\end{array}\right\},
\label{eq:4.9}
\end{equation}

\noindent
where $\nu(\mathbf{A};\pi)=\nu(\mathbf{A})$ is the number of parts of $\pi$ in
color $\mathbf{A}$ and $\lambda(\mathbf{A};\pi)=\lambda(\mathbf{A})$ is the
largest part of $\pi$ in color $\mathbf{A}$, and the other notations in
\eqn{4.9} have similar interpretation.
Now consider partitions $\tilde\pi$ such that part 1 may occur in three
primary colors $\mathbf{A}, \mathbf{B},\mathbf{C}$, but parts $\ge2$ could
occur in the three primary colors as well as in three secondary colors
$\mathbf{AB}, \mathbf{AC}, \mathbf{BC}$ ordered as
\begin{equation}
\mathbf{AB}<\mathbf{AC}<\mathbf{A}<\mathbf{BC}<\mathbf{B}<\mathbf{C}
\label{eq:4.10}
\end{equation}

\noindent
for any given part occurring in these colors and such that the gap
between the parts is $\ge1$ where gap$=1$ only if both parts are either
of the same primary color or if the larger part is in a color of higher
order (as given by \eqn{4.10}).  We call such a partition $\tilde\pi$ a
Type-1 partition as in \cite{AAG}.  It is at this point the interpretation of
the parameters $a,b,c,ab,ac,bc$ becomes clear.  Indeed, we denote
$\nu(\mathbf{A};\tilde\pi)$ by $a$, $\nu(\mathbf{B};\tilde\pi)$ by
$b,\dots,\nu(\mathbf{BC}; \tilde\pi)$ by $bc$.  With this interpretation,
we will now show that $g_{i,j,k}(L,L)$ for $L\ge max(i+j,j+k,k+i)$ is the
generating function for Type-1 partitions $\tilde\pi$ such that
\begin{equation}
\lambda(\tilde\pi)\le L,
\label{eq:4.11}
\end{equation}

\noindent
and the constraints on the frequencies $a,b,\dots,bc$, are as in \eqn{1.2}.
To this end subtract 1 from the smallest part of $\tilde\pi$, 2 from the
second smallest part, $\dots,t$ from the largest part
$\lambda(\tilde\pi)$ of $\tilde\pi$ so that $T_t$ is the total amount
subtracted. (Note that $t=a+b+c+ab+ac+bc$ is the number of parts of
$\tilde\pi$.) Clearly, this subtraction procedure is reversible.
Let the resulting partitions after subtraction be denoted
by $\pi'$.  The colors of the parts of $\pi'$ are those of the parts of
$\tilde\pi$ from which they were derived.  We decompose $\pi'$ into
monochromatic partitions in colors $\mathbf{A}, \mathbf{B}, \mathbf{C},
\mathbf{AB}, \mathbf{AC}, \mathbf{BC}$, denoted as
$\pi'_{\mathbf{A}},\pi'_{\mathbf {B}},\pi'_{\mathbf{C}}, \pi'_{\mathbf{AB}},
\pi'_{\mathbf{AC}},\pi'_{\mathbf{BC}}$.  The monochromatic partitions
satisfy the following conditions:
\begin{equation}
\left.\begin{array}{lll}
\mbox{ 0 could be part of } \pi'_{\mathbf{A}}, & \lambda(\pi'_{\mathbf{A}})=
L-t, &\nu(\pi'_{\mathbf{A}}) = a,\\
\mbox{ 0 could be part of } \pi'_{\mathbf{B}}, & \lambda(\pi'_{\mathbf{B}})=
L-t, &\nu(\pi'_{\mathbf{B}}) = b,\\
\mbox{ 0 could be part of } \pi'_{\mathbf{C}}, & \lambda(\pi'_{\mathbf{C}})=
L-t, &\nu(\pi'_{\mathbf{C}}) = c,\end{array}\right\}
\label{eq:4.12}
\end{equation}
\begin{equation}
\left.\begin{array}{lll} \pi'_{\mathbf{AB}} \mbox{ has distinct parts },&
\lambda(\pi'_{\mathbf{AB}})\le L-t, &\nu(\pi'_{\mathbf{AB}})=ab,\\
\pi'_{\mathbf{AC}} \mbox{ has distinct parts },&
\lambda(\pi'_{\mathbf{AC}})\le L-t, &\nu(\pi'_{\mathbf{AC}})=ac,\\
\pi'_{\mathbf{BC}} \mbox{ has distinct parts },&
\lambda(\pi'_{\mathbf{BC}})\le L-t, &\nu(\pi'_{\mathbf{BC}})=bc,
\end{array}\right\}
\label{eq:4.13}
\end{equation}
and
$$
s(\pi'_{\mathbf{BC}})\ge0, \qquad {\text{if}}\quad s(\pi'_{\mathbf{A}})=0,
\nonumber
$$
$$
s(\pi'_{\mathbf{BC}})\ge1, \qquad {\text{otherwise}},
\nonumber
$$
where $s(\pi'_{\mathbf{BC}})$ is the smallest part in $\pi'_{\mathbf{BC}}$
with similar interpretation for $s(\pi'_{\mathbf{A}})$.
The first summation on the left in \eqn{1.1} is the generating function of
Type--1 partitions $\tilde\pi$ such that $\pi'$ satisfies \eqn{4.12}, \eqn{4.13}
and $s(\pi'_{\mathbf{BC}})\ge1$, $s(\pi'_{\mathbf{A}})\ge0$.
The second summation on the left in \eqn{1.1} is the generating function of
Type--1 partitions $\tilde\pi$ such that $\pi'$ satisfies \eqn{4.12}, \eqn{4.13} and
$s(\pi'_{\mathbf{BC}})=0, s(\pi'_{\mathbf{A}})=0$. Hence, $g_{i,j,k}(L,L)$ is
the generating function of all Type--1 partitions satisfying \eqn{4.11}. Thus,
we have the following new bounded version of the Alladi-Andrews-Gordon \cite{AAG}
refinement of G\"ollnitz's theorem \cite{G}.
\begin{thm}
Let $G_L(n;a,b,c,ab,ac,bc)$ denote the number of Type--1 partitions
$\tilde\pi$ of $n$ such that $\lambda(\tilde\pi)\le L$, $\nu(\mathbf
{A};\tilde\pi)=a,\dots,\nu(\mathbf{BC};\tilde\pi)=bc$.
Let $P_L(n;i,j,k)$ denote the number of partitions $\pi$ of $n$ into
parts occurring in three colors $\mathbf{A}<\mathbf{B}<\mathbf{C}$ such that
parts of the same color are distinct and conditions \eqn{4.9} are
satisfied.  Then for $L\ge max(i+j,j+k,k+i)$ we have
\[
\sum_{i,j,k\mbox{ constraints}} G_L(n;a,b,c,ab,ac,bc)=P_L(n;i,j,k),
\]
where the $i,j,k$ constraints on the summation variables $a,b,c,ab,ac,bc$ are
as in \eqn{1.2}.
\end{thm}

{\bf{Remark 2}}: We would like to stress that the condition \eqn{4.7} is
crucial for our partition theoretical interpretation of $g_{i,j,k}(L,L)$.
Indeed, this condition along with $M=L$ guarantees that whenever the
summands in \eqn{1.1} are non zero, then $L-t \ge 0$. 
Thus the largest parts  $\lambda(\pi'_{\mathbf{A}}),
\dots, \lambda(\pi'_{\mathbf{BC}})$ in \eqn{4.12},  \eqn{4.13} will never take
negative values (see Appendix B for details).

{\bf{Remark 3}}: To see the connection between Theorem 1 and G\"ollnitz's 
theorem \cite{G}, proceed as follows. First, denote part $n$ of color 
${\bf A}$ by $\mathbf{A}_n$, part $n$ of color ${\bf B}$ by 
$\mathbf{B}_n$, etc. Next, replace $\mathbf{A}_n$ by $6n-4$, $\mathbf{B}_n$ 
by $6n-2$, $\mathbf{C}_n$ by $6n-1$, $\mathbf{AB}_n$ by $6n-6$, 
$\mathbf{AC}_n$ by $6n-5$, $\mathbf{BC}_n$ by $6n-3$, let $L \to \infty$, 
and sum over $i,j,k$. This yields G\"ollnitz's (Big) partition theorem 
\cite{G}:

\begin{thm}

Let $B(n)$ denote the number of partitions of $n$ into distinct parts $\equiv$
$2,4,5(mod$ $6)$.\\
Let $C(n)$ denote the number of partitions of $n$ in the
form $m_1+m_2+...+m_s$, no part equals  $1$ or $3$, and such that
$m_l-m_{l+1}\geq 6$ with strict inequality if $m_l \equiv 6,7$ or $9 (mod$
$6)$.

Then, $B(n) = C(n)$.
\end{thm}

\section{Prospects}

In \cite{AAG}, Alladi, Andrews and Gordon discuss companions to G\"ollnitz's
theorem generated by different orderings of the colored integers.  They
show that these companion partition functions are bijectively equivalent
to $G(n;a,b,c,ab,ac,bc)$, and therefore the left hand side of their key
identity \eqn{1.6} is the generating function for all these companion
partition functions.  It turns out that when bounds are imposed on the
parts, these bijections can fail at the boundary.  Thus the finite key
identity \eqn{1.1} in the case $L=M$ corresponds to the bounded G\"ollnitz
partition function in \S4 only with the ordering in \eqn{4.10}.  If a
different ordering were considered as in \cite{AAG}, then this might lead to
a bounded key identity different from \eqn{1.1}, but one which still reduces
to \eqn{1.6} when $L,M\to\infty$.

In a recent paper \cite{AB2}, by studying a reformulation of G\"ollnitz's
theorem as a weighted identity involving partitions into parts differing
by $\ge2$, we deduced several well known results as special cases,
including Jacobi's triple product identity in the form
\begin{equation}
\sum^\infty_{n=-\infty} A^nq^{n^2}=\prod^\infty_{m=1} (1+Aq^{2m-1})
(1+A^{-1} q^{2m-1}) (1-q^{2m}).
\label{eq:5.1}
\end{equation}
Motivated by the method in \cite{AB2} and our Theorem 1 in \S4, we have now
obtained the following new bounded version of \eqn{5.1}
$$
\sum^L_{\ell=0} (-1)^{L+\ell} q^{2(T_L-T_\ell)} \sum^\ell_{n=-\ell}
A^n q^{n^2}=
$$
\begin{equation}
\sum_{\substack{i,j,k\ge0\\ L\ge\max(i+j,i+k,j+k)}} (-1)^k A^{i-j}
q^{2T_i+2T_j+2T_k-i-j}
\left[\begin{array}{c} L-k\\ i\end{array}\right]_{q^2}
\left[\begin{array}{c} L-i\\ j\end{array}\right]_{q^2}
\left[\begin{array}{c} L-j\\ k\end{array}\right]_{q^2}.
\label{eq:5.2}
\end{equation}

\noindent
When $L\to\infty$, \eqn{5.2} reduces to \eqn{5.1}.  The proof and discussion 
of \eqn{5.2} will be presented elsewhere \cite{AB3}. Also in \cite{AB3} we will 
show that \eqn{5.2} implies the new false theta function identity
\begin{equation}
\sum_{\ell\ge0} (-1)^\ell q^{T_\ell}=
\sum_{i,k\ge0}(-1)^{i+k}
\frac{q^{T_i+T_k-ik}}{(q)_i(q)_k}
\left[\begin{array}{c} k+i\\ k\end{array}\right].
\label{eq:5.3}
\end{equation}

Closely related to \eqn{5.1} is another Jacobi's formula
\begin{equation}
\sum_{\ell\ge0}(-1)^{\ell}(2\ell+1) q^{T_\ell}=(q)_{\infty}^3.
\label{eq:5.4}
\end{equation}

\noindent
Using Theorem 1, we found in \cite{AB4} a new polynomial analog of \eqn{5.4} 
\begin{equation}
\sum^L_{\ell=0}(-1)^{\ell} (2\ell+1) q^{T_\ell}=
\sum_{\substack{i,j,k\ge0\\ L\ge\max(i+j,i+k,j+k)}} (-1)^{i+j+k}
q^{T_i+T_j+T_k}
\left[\begin{array}{c} L-k\\ i\end{array}\right]
\left[\begin{array}{c} L-i\\ j\end{array}\right]
\left[\begin{array}{c} L-j\\ k\end{array}\right],
\label{eq:5.5}
\end{equation}

\noindent
which is unexpectedly elegant and is very different from the polynomial 
identity proven by Hirschhorn \cite{H}. Clearly, as $L \to \infty$ 
\eqn{5.5} reduces to \eqn{5.4}. More generally, we derived
\begin{equation}
\sum^L_{\ell=0} a^{-\ell} \frac{1+a^{2\ell+1}}{1+a} q^{T_\ell}=
\sum_{\substack{i,j,k\ge0\\ L\ge\max(i+j,i+k,j+k)}} a^{i-j} (-1)^{k}
q^{T_i+T_j+T_k}
\left[\begin{array}{c} L-k\\ i\end{array}\right]
\left[\begin{array}{c} L-i\\ j\end{array}\right]
\left[\begin{array}{c} L-j\\ k\end{array}\right].
\label{eq:q.carl}
\end{equation}
If we set $q=1$ in \eqn{q.carl} it becomes
\begin{equation}
\frac{a^{L+1}-a^{-L-1}}{a-a^{-1}}=
\sum_{\substack{i,j,k\ge0\\ L\ge\max(i+j,i+k,j+k)}} a^{i-j} (-1)^{k}
\left[\begin{array}{c} L-k\\ i\end{array}\right]_1
\left[\begin{array}{c} L-i\\ j\end{array}\right]_1
\left[\begin{array}{c} L-j\\ k\end{array}\right]_1,
\label{eq:carlitz}
\end{equation}
which is, essentially, a special case of Carlitz's formula (3.7) in \cite{C}
with $n=L$ and $x=y^{-1}=a$.
Actually, from Theorem 1, many new
finite versions of other fundamental results can be deduced, and these
will be presented in \cite{AB3} and \cite{AB4}. Interestingly, all these new
finite identities can be interpreted as $q-$analogs of Carlitz's formulas
for the binomial cycles of length $3$. 

It was the appearance of the
$q-$binomial cycles in \eqn{4.1} in a special form that led to Theorem 1.  In
collaboration with Andrews we intend to conduct a systematic study of
$q-$binomial cycles; in particular we will show that the generating
function of the $q-$binomial cycles of length $3$ can be evaluated in terms
of infinite products.  

Now that we have succeeded in obtaining a partition interpretation of
\eqn{1.1} when $L=M$, it would be worthwhile to see what combinatorial
interpretation \eqn{1.1} has when $L\ne M$. In the case of Schur's partition theorem
with two bounds $L\ne M$ such interpretation turned out to be quite delicate
\cite{AB1}.

Recently, in collaboration with Andrews \cite{AAB}, we have obtained the
following remarkable four parameter key identity:
$$
\sum_{\substack{i,j,k,l\\ constraints}}
\frac{q^{T_t+T_{ab}+T_{ac}+T_{ad}+T_{bc}+T_{bd}+T_{cd}-
bc-bd-cd+4T_{Q-1}+Q(3+2t)}}
{(q)_a(q)_b(q)_c(q)_d(q)_{ab}(q)_{ac}(q)_{ad}(q)_{bc}
(q)_{bd}(q)_{cd}(q)_Q}\times
$$
\begin{equation}
\left\{(1-q^a)+q^{a+bc+bd+Q}(1-q^b)+q^{a+bc+bd+Q+b+cd}\right\}
=\frac{q^{T_i+T_j+T_k+T_\ell}}{(q)_i(q)_j(q)_k(q)_\ell},
\label{eq:5.6}
\end{equation}
where
\begin{equation}
t=a+b+c+d+ab+ac+ad+bc+bd+cd,
\label{eq:5.7}
\end{equation}
and the $i,j,k,l$ constraints on the summation variables $a,b,\dots,cd,Q$ are
\begin{equation}
\left.\begin{array}{l}
i=a+ab+ac+ac+Q,\\
j=b+ab+bc+bd+Q,\\
k=c+ac+bc+cd+Q,\\
\ell=d+ad+bd+cd+Q.\end{array}\right\}
\label{eq:5.8}
\end{equation}

\noindent
Identity \eqn{5.6} reduces to \eqn{1.6} when any one of the parameters
$i,j,k,\ell$ is set equal to 0.  The combinatorial interpretation of
\eqn{5.6} yields a four parameter generalization of G\"ollnitz's theorem.
The discovery and proof of \eqn{5.6} settles a thirty year old problem of
Andrews \cite{A2} who asked whether there exists a partition theorem that lies
``beyond'' the (big) theorem of G\"ollnitz. It would be worthwhile to
seek a bounded identity that reduces to \eqn{5.6} when the bounds go to
infinity, just as \eqn{1.1} reduces to \eqn{1.6} when $L,M\to\infty$.

\paragraph{\large{Note Added}}

Axel Riese informed us that he significantly improved $WZ$ algorithm
and, as a result, was able to obtain a computer proof of the identities
\eqn{1.1} and \eqn{1.6}.

\paragraph{\large{Acknowledgments}}

We wish to thank G. Andrews for a number of stimulating discussions and
for drawing our attention to binomial cycles of Carlitz \cite{C}.
We are grateful to F. Garvan, B. McCoy and A. Riese for their interest and
comments on the manuscript. We would like to thank the referee for his 
comments and suggestions.

\section{Appendix A}

Here we will prove that
\begin{eqnarray}
\left[\begin{array}{c} L\\ s,i,j\end{array}\right]
&=\left[\begin{array}{c} L-1\\
s,i,j\end{array}\right]+q^{L-i}\left[\begin{array}{c} L-1\\
s,i-1,j\end{array}\right]+ q^{L-j} \left[\begin{array}{c} L-1\\
s,i,j-1\end{array}\right]\nonumber\\
&+q^{L-s-i-j} \left[\begin{array}{c} L-1\\
s-1,i,j\end{array}\right]+q^{L-i-j} (1-q^{L-1}) \left[\begin{array}{c}
L-2\\ s,i-1,j-1\end{array}\right].
\label{eq:A.1}
\end{eqnarray}
First we note that
\begin{equation}
\left[\begin{array}{c} L\\
s,i,j\end{array}\right]=\left[\begin{array}{c} L\\
i,j,s\end{array}\right]= \left[\begin{array}{c} L\\
i,j\end{array}\right] \left[\begin{array}{c} L-i-j\\
s\end{array}\right].
\label{eq:A.2}
\end{equation}

\noindent
Next, we recall the symmetric recursion relation
\begin{equation}
\left[\begin{array}{c} L\\ i,j\end{array}\right]=\left[\begin{array}{c}
L-1\\ i,j\end{array}\right]+ q^{L-i}\left[\begin{array}{c} L-1\\
i-1,j\end{array}\right]+ q^{L-j} \left[\begin{array}{c} L-1\\
i,j-1\end{array}\right]+ q^{L-i-j}(1-q^{L-1}) \left[\begin{array}{c}
L-2\\ i-1,j-1\end{array}\right],
\label{eq:A.3}
\end{equation}

\noindent
proven in \cite{AB1}. Combining \eqn{A.2} and \eqn{A.3}, we obtain
\begin{eqnarray}
\left[\begin{array}{c} L\\ s,i,j\end{array}\right]
&=\left[\begin{array}{c} L-1\\ i,j\end{array}\right]
\left[\begin{array}{c} L-i-j\\
s\end{array}\right]+q^{L-i}\left[\begin{array}{c} L-1\\
i-1,j,s\end{array}\right]+q^{L-j} \left[\begin{array}{c} L-1\\
i,j-1,s\end{array}\right]\nonumber\\
&+q^{L-i-j}(1-q^{L-1})\left[\begin{array}{c} L-2\\
i-1,j-1,s\end{array}\right].
\label{eq:A.4}
\end{eqnarray}

\noindent
Finally, using the $q-$binomial recurrence \eqn{2.3} with $m=L-s-i-j$ and
$n=s$, we obtain from \eqn{A.4}
\begin{eqnarray}
\left[\begin{array}{c} L\\ s,i,j\end{array}\right]
&=\left[\begin{array}{c} L-1\\
i,j,s\end{array}\right]+q^{L-s-i-j}\left[\begin{array}{c} L-1\\
i,j,s-1\end{array}\right]+q^{L-i} \left[\begin{array}{c} L-1\\
i-1,j,s\end{array}\right]\nonumber\\
&+q^{L-j}\left[\begin{array}{c} L-1\\
i,j-1,s\end{array}\right]+q^{L-i-j}(1-q^{L-1})\left[\begin{array}{c} L-2\\
i-1,j-1,s\end{array}\right],
\label{eq:A.5}
\end{eqnarray}

\noindent
which is essentially \eqn{A.1}. Formula \eqn{2.11} follows 
from \eqn{A.1} with the substitutions
\begin{equation}
\left.\begin{array}{l} L\mapsto L-s,\\ i\mapsto i-s,\\
j\mapsto j-s.\end{array}\right\}
\label{eq:A.6}
\end{equation}

\section{Appendix B}

Throughout this appendix we assume that $M=L$ in \eqn{1.1}. We need to show
that when \eqn{4.7} holds, $L-t \geq 0$ in all cases where the summands in
\eqn{1.1} are non zero in value. To this end, it is important
to observe that definition \eqn{1.3} implies that 
\begin{eqnarray}
\left[\begin{array}{c} n\\ m \end{array}\right] \neq 0,  \mbox{ iff } m \geq 0
\mbox{ and either } n < 0 \mbox{ or } n \geq m. 
\label{eq:B.1} 
\end{eqnarray}
So, if $n,m \geq 0$, then
\begin{eqnarray}
\left[\begin{array}{c} n\\ m \end{array}\right] \neq 0, 
\mbox{ iff } n \geq m. 
\label{eq:B.2} 
\end{eqnarray}
Let us now assume that $k \leq min(i,j)$ and consider
\begin{eqnarray}
\left[\begin{array}{c} L-t+c\\ c\end{array}\right]
=\left[\begin{array}{c} L-i-j+ab\\ c\end{array}\right],
\label{eq:B.3}
\end{eqnarray}
which appears in the lhs of \eqn{1.1}. Since $L-i-j \geq 0$ (by \eqn{4.7}) and
$ab \geq 0$, then $L-t+c \geq 0$ and, as a result, $L-t \geq 0$ (by \eqn{B.2}).
Obviously, the case $j \leq min(i,k)$ can be treated in the analogous
fashion. If $i \leq min(j,k)$, one needs to consider two $q$-binomials
\begin{eqnarray}
\left[\begin{array}{c} L-t+a-1\\ a-1\end{array}\right]
\left[\begin{array}{c} L-t\\ bc-1\end{array}\right]
=\left[\begin{array}{c} L-j-k+bc-1\\ a-1\end{array}\right]
\left[\begin{array}{c} L-t\\ bc-1\end{array}\right],
\label{eq:B.4}
\end{eqnarray}
which appear in the second product in \eqn{1.1}. Since $L-j-k \geq 0$ (by
\eqn{4.7}) and $bc-1\geq 0$, we conclude that $L-t+a-1 \geq 0$, which again
implies that $L-t \geq 0$.


\begin{thebibliography}{99}

\bibitem{AAB} K. Alladi, G. E. Andrews and A. Berkovich, {\it{ A four
parameter generalization of G\"ollnitz's (Big) partition theorem}},
submitted to Proc. 2000 DIMACS Conf. On Unusual Applications of Number
Theory, (M.B. Nathanson Ed.), CO/0005157.

\bibitem{AAG} K. Alladi, G. E. Andrews, and B. Gordon, {\it{ Generalizations
and refinements of a partition theorem of G\"ollnitz}}, J. Reine Angew
Math, {\bf{460}} (1995) 165-188.

\bibitem{AB1} K. Alladi and A. Berkovich, {\it{ A double bounded version of
Schur's partition theorem}}, to appear in Combinatorica - Erd\"os memorial
issue, CO/0006207.

\bibitem{AB2} K. Alladi and A. Berkovich, {\it{New weighted Rogers-Ramanujan
partition theorems and their implications}}, CO/0009171.

\bibitem{AB3} K. Alladi and A. Berkovich, {\it{ New finite versions of the
Jacobi triple product, Sylvester, and Lebesgue identities}} (in preparation).

\bibitem{AB4} K. Alladi and A. Berkovich, {\it{ New finite
versions of identities of Jacobi and Gauss and of the Hecke modular form}} 
(in preparation).

\bibitem{A1} G. E. Andrews, {\it{On a partition theorem of G\"ollnitz and
related formulae}}, J. Reine Angew. Math., {\bf{236}}(1969), 37-42.

\bibitem{A2} G. E. Andrews, {\it{ The use of computers in the search of
identities of Rogers-Ramanujan type}}, in Computers in Number Theory
(A.O.L. Atkin and B. J. Birch, Eds.) Academic Press (1971), 377-387.

\bibitem{BER} A. Berkovich, {\it{Fermionic counting of RSOS states and Virasoro
character formulas for the unitary minimal series $M(\nu,\nu+1)$: Exact results}},
Nucl. Phys. {\bf{B431}}(1994), 315-348.

\bibitem{BMS} A. Berkovich, B. McCoy, and A. Schilling, {\it{
Rogers-Schur-Ramanujan type identities for the M(p,p') minimal models
of conformal field theory}}, Comm. Math. Phys., {\bf{191}}(1998), 325-395.

\bibitem{C} L. Carlitz, {\it Some multiple sums and binomial identities},
J. Soc. Indust. Appl. Math., {\bf{13}} (1965), 469-486.

\bibitem{GR} G. Gasper and M. Rahman, {\it{ Basic hypergeometric series}},
in Encyclopedia of Mathematics and its Applications, Vol. 10, Cambridge (1990).

\bibitem{G} H. G\"ollnitz, {\it Partitionen mit Differenzenbedingungen}, J.
Reine Angew. Math., {\bf{225}}(1967), 154-190.

\bibitem{H} M. Hirschhorn, {\it Polynomial identities which imply identities
of Euler and Jacobi}, Acta Arith., {\bf{XXXII}} (1977), 73-78.

\bibitem{SW} A. Schilling and S. Ole Warnaar, {\it{Supernomial coefficients,
polynomial identities, and q-series}}, The Ramanujan J., {\bf{2}}(1998),
459-494.

\bibitem{Z} D. Zeilberger, {\it{A q-Foata proof of the q-Saalsch\"utz 
identity}}, Europ. J. Comb., {\bf{8}} (1987), 461-463. 

\end{thebibliography}
\end{document}